\documentclass[11pt]{article}

\usepackage{amsfonts,amsmath,latexsym}
\usepackage{amstext,amssymb}
\usepackage{amsthm}
\usepackage{color}
\usepackage{soul}

\newcommand {\R}{\mathbf{R}}

\newcommand {\B}{\mathbf{B}}

\newcommand {\calCp}{\mathcal{C}_p}

\newtheorem {thm} {Theorem}

\newtheorem {lemma}{Lemma}
\newtheorem {cor} [thm] {Corollary}

\newtheorem {rmk} {Remark}

\newcommand{\beq}{\begin{equation}}
\newcommand{\eeq}{\end{equation}}
\numberwithin{equation}{section}

\begin {document}

\title{ 
{\bf An isoperimetric inequality for extremal Sobolev functions}}
\author{Tom Carroll \footnote{School of Mathematical Sciences, University 
College Cork, {\tt t.carroll@ucc.ie}}
 and Jesse Ratzkin \footnote{Department of Mathematics and Applied 
Mathematics, University of Cape Town, {\tt jesse.ratzkin@uct.ac.za}}
}
\maketitle

\begin {abstract} 
\noindent Let $D \subset \R^n$  be a bounded domain with a Lipschitz boundary, 
let $1< p< \frac{2n}{n-2}$, and let $\phi$ minimize the ratio $\|\nabla u\|_{L^2}/
\|u\|_{L^p}$. We prove a reverse-H\"older inequality, finding a lower bound for 
$\|\phi\|_{L^{p-1}}$ in terms of $\|\phi\|_{L^p}$, in which equality holds if and only if 
$D$ is a ball. This result generalizes an inequality due to Payne and Rayner
\cite{PR,PR2} regarding eigenfunctions of the Laplacian. 
\end {abstract}

\section{Introduction and statement of results}

Let $D \subset \R^n$ be a bounded domain with Lipschitz 
boundary, and let $1 < p <\frac{2n}{n-2}$ (or, $p > 1$ 
if $n=2$). For this range of exponents, the Sobolev 
embedding $W^{1,2}_0(D) \hookrightarrow L^p(D)$ is 
compact, and so the infimum 
\begin {equation}\label{sobolev-quotient} 
\calCp(D) =\inf \left \{ \frac{\int_D |\nabla u |^2 d\mu}
{\left (\int_D |u|^p d\mu\right )^{2/p}}: u \in W^{1,2}_0 (D), 
u \not \equiv 0\right \}\end {equation} 
is finite and achieved by a nontrivial function $\phi = \phi_p$. 

We take 
this opportunity to set notation for the remainder of the paper. We 
denote the volume element of the usual Lebesque measure in 
$\R^n$ by $d\mu$; when it will be necessary, we will denote the induced  area 
element on a hypersurface $\Sigma \subset \R^n$ by $d\sigma$. We write 
the appropriate dimensional volume of a set as $\Omega$ as 
$|\Omega|$, {\it i.e.} if $\Omega \subset \R^n$ is an open set then
$|\Omega| = \mu(\Omega)$ and if $\Sigma\subset \R^n$ is a hypersurface 
then $|\Sigma| = \sigma(\Sigma)$. If 
$\B_1 \subset \R^n$ is the unit ball, we denote $|\B_1| = \omega_n$, 
so that $|\B_r| = \omega_n r^n$ and $|\partial \B_r| = 
n \omega_n r^{n-1}$.  
The Sobolev space $W^{1,2}_0(D)$ is the closure of 
$\mathcal{C}^\infty_0(D)$ under the norm 
$\|u\|_{W^{1,2}} ^2 = \|u\|_{L^2}^2 + \|\nabla u \|_{L^2}^2$.

An extremal function $\phi$ for \eqref{sobolev-quotient} will solve the boundary value problem 
in $D$:
\begin {equation} \label{sobolev-pde} 
\Delta \phi + \lambda\, \phi^{p-1} = 0 , \quad \left. 
\phi \right |_{\partial D} = 0.\end {equation} 
Without loss of generality we can take $\phi>0$  inside $D$. General regularity results imply 
that $\phi \in C^\infty_0(D)$, and a short integration by parts argument reveals that 
\begin {equation} \label{sobolev-scaling} 
\lambda = \calCp(D) \left ( \int_D \phi^p d\mu \right )^{\frac{2-p}{p}}.
\end {equation} 

This sharp Sobolev constant $\calCp(D)$ and its associated extremal 
function $\phi_p$ are both the subject of a vast literature, and incorporate 
much information relating the function theory and the geometry of $D$. 
In particular, a long string of results (for example, \cite{PR,Chiti,Alvino,CFMP})
have uncovered isoperimetric-type inequalities of varying sorts. Our 
main theorem generalizes the reverse-H\"older inequalities 
of \cite{PR,PR2,CRpPR1}, and has the following form. 
\begin {thm} \label{main-thm} 
Let $D \subset \R^n$ be a bounded domain with Lipschitz boundary, 
let $\calCp(D)$ be the sharp Sobolev constant defined by \eqref{sobolev-quotient}, 
and let $\phi$ be its associated extremal function. 
Let $D^*$ be a ball with the same volume as $D$. 
Then 
\begin {equation} \label {main-ineq}
\left ( \int_D \phi^{p-1}\, d\mu \right )^2 \geq  |D|^{\frac{n-2}{n}} 
\left ( \int_D\phi^p\, d\mu \right )^{\frac{2(p-1)}{p}} 
\left [ \frac{2n^2 \omega_n^{2/n}}{p\,\calCp(D)} -
 (n-2)\frac{  n\, \omega_n^{\frac{2}{n} + \frac{p^2-p+2}{p(p-1)}}}{C_p(D^*)}
\right ].
\end {equation}
Equality holds if and only if $D$ is a ball. 
\end {thm} 

\begin {rmk} 
\begin {itemize} 
\item The H\"older inequality implies that for any  
$u \in W^{1,2}_0$ we have 
$$\int_D |u|^{p-1}\, d\mu \leq |D|^{1/p} \left ( \int_D |u|^p \,d\mu \right )^{\frac{p-1}{p}}.$$
For this reason, upper bounds of the form \eqref{main-ineq} are called 
reverse-H\"older inequalities. 
\item Observe that we recover the main inequality of \cite{PR2} in the case 
$p=2$, and we recover the reverse-H\"older inequality of \cite{CRpPR1} in the 
case $n=2$. 
\end {itemize} 
\end {rmk}

\bigskip \noindent {\sc Acknowledgements:} T. C. is partially supported by the programme 
of the ESF Network  \lq Harmonic and Complex Analysis and Applications\rq\ (HCAA).
J. R. is partially supported by a Carnegie Research grant from the University 
of Cape Town. 

\section{Proof of the main theorem} \label{proof-main-thm} 

We begin by briefly outlining our strategy for proving \eqref{main-ineq}, 
which we adapted from Payne and Rayner's   proof in \cite{PR2}. 
Let $M = \sup_{x \in D} \phi(x)$, and for $0 \leq t \leq M$ we 
define 
$$D_t = \{x \in D: \phi (x) > t\}, \qquad \Sigma_t = \{x \in D: \phi (x) = t\}.$$
By Sard's theorem, we have $\Sigma_t = \partial D_t$ for almost every 
value of $t$. To prove \eqref{main-ineq} we define the auxilliary function 
$$H(t) = \int_{D_t} \phi^{p-1} \,d\mu = \int_t^M \tau^{p-1} \int_{\Sigma_\tau}
\frac{d\sigma}{|\nabla \phi|}\, d\tau, \quad t\in [0,M].$$
In Section \ref{diff-ineq} we derive lower bounds for the second derivative 
of $H$, and in Section \ref{int-ineq} we integrate these to obtain 
several integral inequalites for $H$ and for powers of $\phi$. In 
Section \ref{radially-symm} we examine the one-dimensional eigenvalue problem which 
arises from the radially symmetric case, and in Section \ref{main proof} complete the proof of \eqref{main-ineq}.

\subsection{Differential inequalities} \label{diff-ineq} 
We let $V(t) = |D_t|$. Then, by the co-area formula, 
$$V'(t) = -\int_{\Sigma_t}\frac{d\sigma}{|\nabla \phi|} < 0.$$
 Thus $V$ is a monotone function of $t$, and we can invert it to 
obtain $t = t(V)$, with 
$$\frac{dt}{dV} = \frac{1}{V'(t)} = -\frac{1}{\int_{\Sigma_t} \frac{d\sigma}
{|\nabla \phi|}}.$$ 
This in turn implies that
\begin {equation}\label{change-var1}
\frac{dH}{dV} = \frac{dH}{dt}\frac{dt}{dV} = \left ( -t^{p-1} \int_{\Sigma_t} \frac{d\sigma}
{|\nabla \phi|}\right ) \cdot\left (- \frac{1}{\int_{\Sigma_t} \frac{d\sigma}{|\nabla \phi|}}
\right ) = t^{p-1}, \end {equation} 
a relation which will prove quite useful in our computations. 
Taking one more derivative shows that 
$$\frac{d^2H}{dV^2} = \frac{d}{dV} \big(t^{p-1}\big) = -\frac{(p-1)\,t^{p-2}}
{\int_{\Sigma_t} \frac{d\sigma}{|\nabla \phi|}}.$$ 

\begin {lemma} 
The function $H$ satisfies 
\begin {equation} \label {comparison-ode1} 
\frac{d^2 H}{dV^2}  \geq -(p-1) (t(V))^{p-2} \frac{\Lambda H(V)}
{n^2 \omega_n^{2/n} V^{\frac{2(n-1)}{n}}},
\quad V \in \big[ 0,\vert D \vert \big].
\end {equation} 
with the boundary conditions $H(0) = 0$ and $H'(|D|) = 0$. Moreover, 
equality in \eqref{comparison-ode1} forces $D$ to be a ball, and 
forces the function $\phi$ to be radially symmetric. 
\end {lemma} 

\begin {proof} 
By the Cauchy-Schwarz inequality,
$$|\Sigma_t|^2 \leq \left ( \int_{\Sigma_t} |\nabla \phi| \,d\sigma \right ) 
\left (\int_{\Sigma_t} \frac{d\sigma}{|\nabla \phi|} \right ),$$
which we can rearrange to read 
\begin {equation} \label {cauchy-schwarz1}
\int_{\Sigma_t} \frac{d\sigma}{|\nabla \phi |} \ \geq\ \frac{|\Sigma_t|^2}
{\int_{\Sigma_t} |\nabla \phi| \,d\sigma}.\end {equation} 
Since $\Sigma_t$ is a level-set of $\phi$,  we may use 
the divergence theorem and \eqref{sobolev-pde} to obtain 
\begin {eqnarray*}
\int_{\Sigma_t} |\nabla \phi|\,d\sigma & = & -\int_{\Sigma_t} \frac{\partial \phi}
{\partial \eta} \,d\sigma = - \int_{D_t} \Delta \phi \,d\mu \\ 
& = & \lambda \int_{D_t} \phi^{p-1}\,d\mu = \lambda H(t).
\end {eqnarray*} 
Combining this with \eqref{cauchy-schwarz1} we obtain  
\begin {equation}\label{cauchy-schwarz2}
\int_{\Sigma_t} \frac{d\sigma}{|\nabla \phi |} \ \geq\ \frac{|\Sigma_t|^2}
{\lambda\,H(t)}.\end{equation} 
By the classical isoperimetric inequality, 
\begin {equation}\label{isop-ineq1}
|\Sigma_t|^2 \geq n^2 \omega_n^{2/n} |D_t|^{\frac{2(n-1)}{n}}
= n^2 \omega_n^{2/n} V^{\frac{2(n-1)}{n}}.
\end {equation}
Together with \eqref{cauchy-schwarz2}, this shows that
\begin {eqnarray*} 
\frac{d^2 H}{dV^2} & = & - (p-1) \frac{t^{p-2}}{\int_{\Sigma_t}
\frac{d\sigma}{|\nabla \phi|}} \geq -(p-1)\, t^{p-2} \,\frac{\lambda H(V)}
{|\Sigma_t|^2} \\ 
& \geq & -(p-1)\,t^{p-2} \,\frac{\lambda H(V)}{n^2 \omega_n^{2/n}
V^\frac{2(n-1)}{n}}.
 \end {eqnarray*} 
Notice that the boundary conditions for this differential inequality 
are 
\begin{equation}\label{H-bc-V}
H(0) = 0, \qquad H' (|D|) =\left.  t^{p-1} \right |_{t=0} = 0.
\end{equation}
Moreover, we only have equality in \eqref{comparison-ode1}
for each $V$ in $\big[ 0,\vert D \vert \big]$ if 
we have equality in \eqref{isop-ineq1} for almost every $t$, which 
in turn implies that $\Sigma_t$ is a round sphere for almost 
every $t \in [0,M]$. This is possible only if 
$D$ is itself a ball. Also, equality in \eqref{comparison-ode1} forces 
equality in \eqref{cauchy-schwarz1}, which implies 
$\nabla \phi$ must be constant on each sphere $\Sigma_t$, and so 
$\phi$ must be radial. \end {proof} 

We change variables by letting $\rho = (V/\omega_n)^{1/n}$ 
be the volume radius of $D_t$, so that $V = |D_t| = \omega_n \rho^n$.  
We also define $\rho_M = (|D|/\omega_n)^{1/n}$. As a function of $\rho$, the function 
$H$ satisfies the boundary conditions
\begin{equation}\label{H-bc-rho}
H(0) = H'(0) = \cdots = H^{(n-1)}(0) = 0, \qquad H'(\rho_M) = 0.
\end{equation}

\begin {lemma} 
\begin {equation} \label{comparison-ode2}
\frac{d}{d\rho} \left [\rho^{1-n} \left( \frac{dH}{d\rho} \right )^{\frac{1}{p-1}} \right ] 
\geq -\frac{\lambda}{(n\omega_n)^{\frac{p-2}{p-1}}} \,\rho^{1-n} \,H(\rho), 
\quad 0 < \rho < \rho_M.
\end {equation} 
\end {lemma} 
\begin {proof} 
Taking derivatives, we see that 
\begin {equation} \label{change-var2} 
\frac{dH}{d V} = \frac{\rho^{1-n}}{n\omega_n} \frac{d H}{d \rho}, 
\quad \frac{d^2 H }{d\rho^2} = \frac{\rho^{1-n}}{n\omega_n}
\frac{d}{d\rho} \left (\frac{\rho^{1-n}}{n\omega_n} \frac{dH}{d\rho}
\right ).\end {equation} 
Substituting  these expressions in \eqref{comparison-ode1} gives 
\begin {equation} \label{comparison-ode3}
\frac{d}{d\rho} \left ( \rho^{1-n} \frac{dH}{d\rho} \right )
\geq -(p-1)\,\lambda\, t^{p-2}\, \rho^{1-n}\, H(\rho).
\end {equation} 
However, 
\[
t^{p-1} = \frac{dH}{dV} = \frac{\rho^{1-n}}{n\omega_n} \,\frac{dH}{d\rho},
\]
so that \eqref{comparison-ode3} becomes 
\[
\frac{d}{d\rho} \left ( \rho^{1-n} \,\frac{dH}{d\rho} \right )\geq 
-\frac{(p-1)\, \lambda}{(n\omega_n)^{\frac{p-2}{p-1}}} \left ( \rho^{1-n}
\,\frac{dH}{d\rho} \right )^{\frac{p-2}{p-1}} \rho^{1-n}\,H(\rho).
\]
This we can rewrite as 
\[
\frac{1}{p-1} \frac{\frac{d}{d\rho} \left ( \rho^{1-n} \frac{dH}{d\rho} \right )}
{\left ( \rho^{1-n} \frac{dH}{d\rho} \right )^{\frac{p-2}{p-1}} } 
= \frac{d}{d\rho} \left [ \rho^{1-n} \left (\frac{dH}{d\rho} \right )^{\frac{1}{p-1}}
\right ] \geq -\frac{\lambda}{(n\omega_n)^{\frac{p-2}{p-1}}}\, \rho^{1-n}\, H(\rho).
\qedhere
\] 
\end {proof} 

\begin {rmk} Since \eqref{comparison-ode2} is really the 
same as \eqref{comparison-ode1} rewritten in different variables, 
equality holds in \eqref{comparison-ode2} for $0 < \rho < \rho_M$ 
if and only if $D$ is a ball and $\phi$ is radial. 
\end {rmk} 

\subsection{Integral inequalities} \label{int-ineq} 

In this section we integrate \eqref{comparison-ode1} and 
\eqref{comparison-ode2} to obtain inequalities for the integral 
of $H$ and the integral of powers of $\phi$. As each of these 
inequalities is an integrated form of \eqref{comparison-ode1} 
and \eqref{comparison-ode2}, equality holds if and only 
if $D$ is a ball and $\phi$ is radial. 

\begin {lemma} 
\begin {eqnarray} \label {comparison-ode4} 
\left ( \int_D \phi^{p-1} \,d\mu \right )^2 &\geq &\frac{2n^2 \omega_n^{2/n}
|D|^{\frac{n-2}{n}}}{p\, \calCp(D)} \left ( \int_D \phi^p\,d\mu \right)^{\frac{2(p-1)}{p} }\\ \nonumber 
&&\qquad\qquad - \frac{n-2}{n} |D|^{\frac{n-2}{n}} \int_0^{|D|}
V^{\frac{2(1-n)}{n}} H^2(V) \,dV. 
\end {eqnarray}
\end {lemma} 

\begin {proof} 
We multiply the inequality \eqref{comparison-ode1} by 
$\frac{p}{p-1} V \left ( \frac{dH}{dV} \right )^{1/(p-1)}$ 
and integrate from $0$ to $|D|$. Upon integration, the left hand side becomes 
\begin {eqnarray*} 
\int_0^{|D|} \frac{p}{p-1} V \left ( \frac{dH}{dV} 
\right )^{1/(p-1)} \frac{d^2 H}{dV^2}\,dV & = & \int_0^{|D|}
V \frac{d}{dV} \left[\left (\frac{dH}{dV}\right ) ^{p/(p-1)}\right]\, dV \\ 
& = & V \left. \left ( \frac{dH}{dV} \right )^{p/(p-1)} \right|_0^{|D|}
- \int_0^{|D|} \left (\frac{dH}{dV} \right )^{p/(p-1)} dV \\  
& = & - \int_0^{|D|} \big(t^{p-1} (V)\big)^{p/(p-1)} dV \\  
& = & -\int_0^{|D|} t^p(V) dV = -\int_D\phi^p d\mu.
\end {eqnarray*} 
The boundary terms in the integration by parts vanished since $H'(\vert D \vert)=0$, 
while (2.1) 
was used at the third step.
On the other hand, using \eqref{change-var1} again, the right hand side becomes
\begin{align*} 
-\frac{p\,\lambda}{n^2 \omega_n^{n/2}} \int_0^{|D|}
V \left (\frac{dH}{dV}\right ) ^{1/(p-1)} & t^{p-2}  H(t) V^{\frac{2(1-n)}{n}}\, dV\\
& =  -\frac{p\,\lambda}{n^2 \omega_n^{n/2}} \int_0^{|D|} t^{p-2}
\left ( \frac{dH}{dV} \right )^{\frac{2-p}{p-1}} V^{\frac{2-n}{n}} H(V) 
\frac{dH}{dV}\, dV \\ 
& =  -\frac{p\,\lambda}{n^2 \omega_n^{n/2}} \int_0^{|D|}t^{p-2}
(t^{p-1})^{\frac{2-p}{p-1}} V^{\frac{2-n}{n}} H(V) \frac{dH}{dV}\, dV \\ 
& =  -\frac{p\,\lambda}{n^2 \omega_n^{n/2}}\int_0^{|D|} V^{\frac{2-n}{n}}
H(V) \frac{dH}{dV}\, dV
\end{align*}
We combine these last two equations and replace $\lambda$ 
by $\calCp(D) \left (\int_D \phi^p d\mu \right ) ^{(2-p)/p}$ to obtain
\begin {eqnarray*}
\left ( \int_D \phi^p \,d\mu \right )^{\frac{2(p-1)}{p}} & \leq & 
\frac{p\,\calCp(D)}{n^2 \omega_n^{n/2}}
\int_0^{|D|} V^{\frac{2-n}{n}} H(V) \,\frac{dH}{dV}\, dV\\
&  = & \frac{p\, \calCp(D)}{2n^2 \omega_n^{2/n}}\int_0^{|D|}
V^{\frac{2-n}{n}} \frac{d}{dV} \big(H^2(V)\big) \,dV \\ 
& = & \frac{p\,\calCp(D)}{2n^2 \omega_n^{2/n}} \left [ 
|D|^{\frac{2-n}{n}} \left ( \int_D \phi^{p-1} d\mu \right )^2 + 
\frac{n-2}{n} \int_0^{|D|} H^2(V) \,V^{\frac{2(1-n)}{n}}\, dV \right ],
\end {eqnarray*}
which we can rearrange to give \eqref{comparison-ode4}.
\end{proof} 

\begin {lemma} 
\begin {equation} \label {comparison-ode5}
\int_0^{\rho_M} \rho^{\frac{1-n}{p-1}} \left ( \frac{dH}{d\rho}
\right )^{\frac{p}{p-1}} \leq \frac{\lambda}{(n\omega_n
)^{\frac{p-2}{p-1}}} \int_0^{\rho_M} \rho^{1-n} H^2(\rho) d\rho.
\end {equation}
Equality holds if and only if $D$ is a ball. 
\end {lemma} 

\begin {proof}
We mutliply \eqref{comparison-ode2} by $H$ and integrate 
from $0$ to $\rho_M$. The boundary conditions \eqref{H-bc-rho} imply that 
$\rho^{1-n}\,\frac{dH}{d\rho}$ is bounded at $0$. 
Hence the boundary terms vanish in the integration parts below, and we obtain that  
\[
\int_0^{\rho_M}\rho^{\frac{1-n}{p-1}} \left ( \frac{dH}{d\rho} \right )^{\frac{p}{p-1}}\,d\rho 
=
\int_0^{\rho_M} \left [ \rho^{1-n}\,\frac{dH}{d\rho} \right]^{\frac{1}{p-1}} \frac{dH}{d\rho}\, d\rho 
\leq \frac{\lambda}{(n\omega_n)^{\frac{p-2}{p-1}}} \int_0^{\rho_M} \rho^{1-n} H^2(\rho) d\rho.
\]
\end {proof} 

\begin {lemma} 
With $\phi$, $H$, and $\rho$ defined as above, 
\begin {equation} \label {change-var3} 
\int_0^{\rho_M} \rho^{\frac{1-n}{p-1}} \left ( \frac{dH}{d\rho}
\right )^{\frac{p}{p-1}} \,d\rho= (n\omega_n)^{\frac{1}{p-1}}\int_D 
\phi^p\, d\mu.\end{equation} 
\end{lemma} 

\begin {proof} We use \eqref{change-var1} and \eqref{change-var2} 
to conclude 
\begin {eqnarray*} 
\int_0^{\rho_M} \rho^{\frac{1-n}{p-1}} \left ( \frac{dH}
{d\rho} \right )^{\frac{p}{p-1}} & = & \int_0^{\rho_M}
(\rho^{1-n})^{\frac{p}{p-1}} \left ( \frac{dH}{d\rho} \right 
)^{\frac{p}{p-1}} \rho^{n-1} d\rho \\ 
& = & \int_0^{|D|}\left (\rho^{1-n} \frac{dH}{d\rho} 
\right )^{\frac{p}{p-1}} \frac{dV}{n\omega_n}  \\
&= &\int_0^{|D|} \left ( n\omega_n\frac{dH}{dV} \right 
)^{\frac{p}{p-1}} \frac{dV}{n\omega_n} \\ 
& = & (n\omega_n)^{\frac{1}{p-1}} \int_0^{|D|}
\left (\frac{dH}{dV} \right )^{\frac{p}{p-1}}dV\\
& =& (n\omega_n)^{\frac{1}{p-1}}\int_0^{|D|} 
(t^{p-1})^{\frac{p}{p-1}} dV \\ 
& = & (n\omega_n)^{\frac{1}{p-1}} \int_0^{|D|}t^p dV
= (n\omega_n)^{\frac{1}{p-1}} \int_D \phi^p d\mu.
\qedhere
\end {eqnarray*}
\end {proof}

\begin {cor} 
\begin {equation} \label{change-var4} 
n\omega_n \left ( \int_D \phi^p \,d\mu \right )^{\frac{2(p-1)}{p}}\leq
\ \calCp(D)\, \int_0^{\rho_M} \rho^{1-n} H^2(\rho) \,d\rho.
\end {equation}
Moreover, we have equality if and only if $D$ is a ball and $\phi$ is radial. 
\end {cor} 
\begin {proof} Combine \eqref{comparison-ode5}, \eqref{change-var3}, 
and \eqref{sobolev-scaling}. 
\end {proof} 

\begin{lemma} 
\begin {eqnarray} \label {comparison-ode6} 
\left ( \int_D \phi^{p-1} \,d\mu \right )^2 & \geq & 
\frac{2n^2 \omega_n^{2/n} |D|^{\frac{n-2}{n}}}
{p \, \calCp(D)} \left (\int_D \phi^p d\mu \right )^{\frac{2(p-1)}{p}}
\\ \nonumber 
&&\qquad\qquad- (n-2)\,\omega_n^{\frac{2-n}{n}} |D|^{\frac{n-2}{n}}
\int_0^{\rho_M} \rho^{1-n} H^2(\rho) \,d\rho.
\end {eqnarray}
Equality holds if and only if $D$ is a ball.
\end {lemma} 

\begin {proof} 
Since $\rho(V) = \big( V/\omega_n\big)^{1/n}$, we have 
\[
n \omega_n^{1/n} \frac{d \rho}{dV} V^{\frac{n-1}{n}} = 1,
\]
so that 
\begin {eqnarray*} 
\int_0^{|D|} V^{\frac{2(1-n)}{n}} H^2(V)  \,dV& = & 
\int_0^{|D|} H^2(\rho) \,\omega_n^{\frac{2(1-n)}{n}} \rho^{2(1-n)}
\,n \omega_n^{1/n}\, (\omega_n \rho^n)^{\frac{n-1}{n}} \frac{d\rho}{dV} \,dV \\ 
& = & n \omega_n^{\frac{2-n}{n}} \int_0^{\rho_M} \rho^{1-n} H^2(\rho) \,d\rho .
\end {eqnarray*}
Putting this into \eqref{comparison-ode4} gives \eqref{comparison-ode6}.
\end {proof}

\subsection {The radially symmetric case}
\label {radially-symm}

Motivated by \eqref{comparison-ode5} and \eqref{H-bc-rho}, we define $\Lambda_*$ by 
\begin {equation} \label{lambda-star}
\Lambda_* = \inf \left \{ \left. \left (\int_0^{\rho_M} \rho^{\frac{1-n}{p-1}} 
\,f'(\rho)^{\frac{p}{p-1}} \,d\rho \right )^{\frac{2(p-1)}{p}} \right/
\int_0^{\rho_M} \rho^{1-n}\, f^2 (\rho)\, d\rho \right \}
\end{equation}
where the infimum is over all functions on $[0,\rho_M]$ for which
\begin{equation}\label{lambda-star-bc} 
f(0) = f'(0) = \cdots = f^{(n-1)}(0)=0 = f'(\rho_M),\quad f \not \equiv 0 . 
\end {equation} 
%
%
\begin {rmk} Notice that we have rescaled the numerator 
to make the quotient scale-invariant. This does not, however, 
affect the Euler-Lagrange equation involved. 
\end {rmk}

\begin {lemma} The Euler-Lagrange equation for the variational problem 
\eqref{lambda-star}, with the boundary conditions \eqref{lambda-star-bc}, is 
\begin {equation} \label{Euler-Lagrange}
f''(\rho) - \frac{n-1}{\rho}\,f'(\rho) + \Lambda \big[ \rho^{1-n} 
f'(\rho) \big]^{\frac{p-2}{p-1}}\, f(\rho) = 0.\end {equation} 
\end {lemma} 

\begin {proof} Since the ratio defining $\Lambda_*$ is scale-invariant,  we 
may either restrict our attention to either of the constrained critical point problems: 
\[
\text{minimize }\int_0^{\rho_M} \rho^{\frac{1-n}{p-1}} 
\,f'(\rho)^{\frac{p}{p-1}}\,d\rho\ \text{ subject to } 
\int_0^{\rho_M} \rho^{1-n} f^2 d\rho\ = \text{ constant}
\] 
or 
\[
\text{maximize } \int_0^{\rho_M}\rho^{1-n} f^2 d\rho\ 
\text{ subject to }\int_0^{\rho_M} \rho^{\frac{1-n}{p-1}} 
\,f'(\rho)^{\frac{p}{p-1}}\,d\rho\ = \text{ constant.}
\]
Regardless, the method of Lagrange multipliers implies that 
a constrained critical point $f$ satisfies 
\[
\left. \frac{d}{d\epsilon} \right |_{\epsilon = 0} \int_0^{\rho_M}
\rho^{\frac{1-n}{p-1}} \left ( \frac{df}{d\rho} + \epsilon \frac{dg}{d\rho}
\right )^{\frac{p}{p-1}}d\rho = \Lambda \left. \frac{d}{d\epsilon} 
\right |_{\epsilon = 0} \int_0^{\rho_M} \rho^{1-n}\,\big[f(\rho)+ \epsilon\, g(\rho)\big]^2 
d\rho,
\]
for any admissible $g$. 
Having evaluated these derivatives, we use the boundary conditions \eqref{lambda-star-bc} 
to see that $\rho^{1-n}\,f'(\rho)$ is bounded at $0$ and that consequently the boundary terms 
arising from integration by parts vanish, and see that
\begin {align*} 
2\lambda \int_0^{\rho_M} & \rho^{1-n}\, f(\rho)\,g(\rho)\, d\rho \\
&=  
\frac{p}{p-1} \int_0^{\rho_M} \rho^{\frac{1-n}{p-1}} 
\left ( \frac{df}{d\rho} \right )^{\frac{1}{p-1}} \frac{dg}{d\rho}\, d\rho \\ 
& =  -\frac{p}{p-1} \int_0^{\rho_M} g(\rho) \frac{d}{d\rho} \left [ 
\rho^{\frac{1-n}{p-1}} \left (\frac{df}{d\rho} \right )^{\frac{1}{p-1}} 
\right ] \, d \rho \\
& =  -\frac{p}{p-1} \int_0^{\rho_M} g(\rho) \left [ \frac{1}{p-1} 
\,\rho^{\frac{1-n}{p-1}} \left ( \frac{df}{d\rho} \right )^{\frac{2-p}{p-1}}
\frac{d^2 f}{d\rho^2} + \frac{1-n}{p-1}\, \rho^{\frac{2-p-n}{p-1}}
\left ( \frac{df}{d\rho} \right )^{\frac{1}{p-1}} \right ]\,d\rho.
\end {align*}
This must hold for all choices of $g$,  hence (absorbing a factor of 
$2p/(p-1)^2$ into the Lagrange multiplier $\Lambda$) we must have 
\begin {eqnarray*} 
0 & = & \rho^{\frac{1-n}{p-1}} f'(\rho)^{\frac{2-p}{p-1}} f''(\rho) - (n-1)\, 
\rho^{\frac{2-p-n}{p-1}}f'(\rho)^{\frac{1}{p-1}} + \Lambda \, \rho^{1-n}\,f(\rho)\\ 
& = & \rho^{\frac{1-n}{p-1}} f'(\rho)^{\frac{2-p}{p-1}} 
\left [ f''(\rho) - (n-1) \rho^{-1} f'(\rho) + \Lambda \left [ \rho^{1-n}\, f'(\rho) \right ]^{\frac{p-2}{p-1}}
f(\rho) \right ],
\end {eqnarray*}
as claimed.
\end {proof} 

\begin {lemma} 
Let $D^*$ be the ball $\B_{\rho_M}$ of radius $\rho_M$. Then, 
\begin {equation} \label {comparison-ode7} 
\Lambda_* \leq (n \omega_n)^{\frac{2-p}{p}} \calCp(D^*).
\end {equation} 
\end {lemma} 

\begin {proof} We use the function $H(\rho)$ for the ball $\B_{\rho_M}$ as a test function for the 
quotient defining $\Lambda_*$ and use the inequalities  
\eqref{comparison-ode5}, \eqref{sobolev-scaling}, and 
\eqref{change-var3}: 
\begin {align*} 
\Lambda_* & \leq  \left.\left ( \int_0^{\rho_M} \rho^{\frac{1-n}{p-1}}
H'(\rho)^{\frac{p}{p-1}} \,d\rho \right )^{\frac{2(p-1)}{p}} \right/ 
\int_0^{\rho_M} \rho^{1-n}\, H^2(\rho) \,d\rho \\ 
& \leq  \frac{\Lambda}{(n\omega_n)^{\frac{p-2}{p-1}}}
\left ( \int_0^{\rho_M} \rho^{\frac{1-n}{p-1}} H'(\rho)^{\frac{p}{p-1}}\, d\rho \right )^{\frac{p-2}{p}} \\ 
& =  \frac{1}{(n\omega_n)^{\frac{p-2}{p-1}}}\, \calCp(D^*) 
\left ( \int_{D^*} \phi^p\, d\mu  \right )^{\frac{2-p}{p}}
\left [ \frac{1}{(n\omega_n)^{\frac{1}{p-1}}}  \int_{D^*}
\phi^p \,d\mu \right ]^{\frac{p-2}{p}} \\ 
& =  (n\omega_n)^{\frac{2-p}{p}}\, \calCp(D^*). \qedhere
\end {align*} 
\end {proof} 

In order to obtain a lower bound for $\Lambda_*$ in terms of  $\calCp(D^*)$,
we first need to relate the particular $\Lambda$ occurring in the Euler-Lagrage equation \eqref{Euler-Lagrange} 
to the eigenvalue $\Lambda_*$, just as \eqref{sobolev-scaling} relates the number $\lambda$ occurring in 
the Euler-Lagrange equation \eqref{sobolev-pde} to the eigenvalue $\calCp(D)$.

\begin{lemma}\label{star-lemma} 
Let $f$ be a minimizer for $\Lambda_*$ given by \eqref{lambda-star} with the boundary conditions 
\eqref{lambda-star-bc} and satisfy the Euler-Lagrange equation \eqref{Euler-Lagrange}, written as 
\begin{equation}\label{Euler-Lagrange-2}
\frac{d}{d\rho} \left[ \rho^{\frac{1-n}{p-1}} f'(\rho)^{\frac{1}{p-1}}\right] + \Lambda\, \rho^{1-n}\,f(\rho) = 0.
\end{equation}
Then
\begin{equation}\label{star-scaling}
\Lambda = \Lambda_* \, \left (\int_0^{\rho_M} \rho^{\frac{1-n}{p-1}} 
\,f'(\rho)^{\frac{p}{p-1}} \,d\rho \right )^{\frac{2-p}{p}}.
\end{equation}
\end{lemma}
\begin{proof}Multiply the Euler-Lagrange equation\eqref{Euler-Lagrange-2} across by $f(\rho)$ and 
integrate from $0$ to $\rho_M$ to obtain
\[
\int_0^{\rho_M}f(\rho)\,\frac{d}{d\rho} \left[ \rho^{\frac{1-n}{p-1}} f'(\rho)^{\frac{1}{p-1}}\right]\,d\rho
 + \Lambda\, \int_0^{\rho_M}\rho^{1-n}\,f(\rho)^2\,d\rho = 0.
\]
Integrating by parts in the first term and using the boundary conditions \eqref{lambda-star-bc} gives
\[
\int_0^{\rho_M}f(\rho)\,\frac{d}{d\rho} \left[ \rho^{\frac{1-n}{p-1}} f'(\rho)^{\frac{1}{p-1}}\right]\,d\rho
=
- \int_0^{\rho_M} \rho^{\frac{1-n}{p-1}} f'(\rho)^{\frac{p}{p-1}}\,d\rho,
\]
from which it follows that
\[
\Lambda = \left. \int_0^{\rho_M} \rho^{\frac{1-n}{p-1}} f'(\rho)^{\frac{p}{p-1}}\,d\rho 
\right/ \int_0^{\rho_M}\rho^{1-n}\,f(\rho)^2\,d\rho.
\]
We can use \eqref{lambda-star} to write $\int_0^{\rho_M}\rho^{1-n}\,f^2(\rho)\,d\rho$ in terms of $\Lambda_*$ 
since $f$ is a minimizer for this Rayleigh quotient, leading to 
\[
\Lambda = \Lambda_* 
\left( \int_0^{\rho_M} \rho^{\frac{1-n}{p-1}} f'(\rho)^{\frac{p}{p-1}}\,d\rho \right)^{1-\frac{2(p-1)}{p}},
\]
which is \eqref{star-scaling}.
\end{proof}

\begin {lemma} 
\begin {equation} \label {lambda-star2}
\calCp(D^*) \leq  (n \omega_n)^{\frac{p-2}{p}} \Lambda_*.
\end {equation} 
\end {lemma} 

\begin {proof} 
Let $f$ be a minimizer for the generalized quotient \eqref{lambda-star} defining $\Lambda_*$. 
Set 
\[
\psi(\rho) = \int_\rho^{\rho_M} r^{1-n}\, f(r)\,dr, \qquad 0 \leq \rho \leq \rho_M,
\]
so that $\psi(\rho_M) = 0$. Then $\psi(\rho)$ (where $\rho = \vert x \vert$ for $x \in \calCp(D^*)$) 
is an admissible test function for the quotient defining $\calCp(D^*)$. Thus
\begin{equation}\label{t1}
\calCp(D^*) \leq (n \omega_n)^{\frac{p-2}{p}} 
\left. \int_0^{\rho_M}\rho^{n-1}\,\psi'(\rho)^2\,d\rho \right/
\left(\int_0^{\rho_M} \rho^{n-1}\,\psi(\rho)^p\,d\rho\right)^{2/p}.
\end{equation}
Now
\begin{equation}\label{t2}
\int_0^{\rho_M}\rho^{n-1}\,\psi'(\rho)^2\,d\rho 
 = \int_0^{\rho_M}\rho^{n-1}\,\left[ \rho^{1-n}\, f(\rho) \right]^2\,d\rho 
 = \int_0^{\rho_M} \rho^{1-n}\, f(\rho)^2\,d\rho.
\end{equation}
Next, using the Euler-Lagrange equation \eqref{Euler-Lagrange-2}, 
\begin{align*}
\psi(\rho) & = \int_\rho^{\rho_M} r^{1-n}\, f(r)\,dr \\
& = -\frac{1}{\Lambda} \int_\rho^{\rho_M} 
\frac{d}{dr} \left[ r^{\frac{1-n}{p-1}} f'(r)^{\frac{1}{p-1}}\right]\,dr \\
& = -\frac{1}{\Lambda}\,\left. r^{\frac{1-n}{p-1}}\, f'(r)^{\frac{1}{p-1}}\right\vert_{r=\rho}^{r=\rho_M}\\
& = \frac{1}{\Lambda}\, \rho^{\frac{1-n}{p-1}}\, f'(\rho)^{\frac{1}{p-1}},
\end{align*}
where we used $f'(\rho_M)=0$. From this we obtain that 
\begin{align}
\int_0^{\rho_M} \rho^{n-1}\,\psi(\rho)^p\,d\rho 
& = \int_0^{\rho_M} \rho^{n-1}\, \frac{1}{\Lambda^p}\, \rho^{\frac{p(1-n)}{p-1}} f'(\rho)^{\frac{p}{p-1}}\,d\rho
\nonumber \\
& = \frac{1}{\Lambda^p}\, \int_0^{\rho_M} \rho^{\frac{1-n}{p-1}}\, f'(\rho)^{\frac{p}{p-1}}\,d\rho.
\label{t3}
\end{align}
With the help of the identities \eqref{t2} and \eqref{t3}, we can write the numerator and the denominator
of \eqref{t1} in terms of the minimizer $f$ for $\Lambda_*$. 
We find, using that $f$ minimizes the quotient for $\Lambda_*$ at the second step and using 
\eqref{star-scaling} at the third step, that 
\begin{align}
\calCp(D^*) & \leq (n \omega_n)^{\frac{p-2}{p}}  
	\left. \int_0^{\rho_M} \rho^{1-n}\, f(\rho)^2\,d\rho\right/ 
	\left(\frac{1}{\Lambda^p}\, \int_0^{\rho_M} \rho^{\frac{1-n}{p-1}}\, f'(\rho)^{\frac{p}{p-1}}\,d\rho\right)^{\frac{2}{p}}\\
& = (n \omega_n)^{\frac{p-2}{p}}\,  \frac{\Lambda^2}{\Lambda_*}\, 
\left (\int_0^{\rho_M} \rho^{\frac{1-n}{p-1}} 
\,f'(\rho)^{\frac{p}{p-1}} \,d\rho \right )^{\frac{2(p-1)}{p}-\frac{2}{p}} \\
& = (n \omega_n)^{\frac{p-2}{p}}\,  \frac{1}{\Lambda_*}\, 
\left[ \Lambda\, \left (\int_0^{\rho_M} \rho^{\frac{1-n}{p-1}} 
\,f'(\rho)^{\frac{p}{p-1}} \,d\rho \right )^{\frac{p-2}{p}}\right]^2\\
& = (n \omega_n)^{\frac{p-2}{p}}\,  \frac{\Lambda_*^2}{\Lambda_*} 
=  (n \omega_n)^{\frac{p-2}{p}}\, \Lambda_*.
\qedhere
\end{align}
\end{proof}

\subsection{Completion of the proof of Theorem~\ref{main-thm}}
\label{main proof}
We can now finally complete the proof of Theorem \ref{main-thm}. 
Indeed, we have 
\begin {eqnarray*} 
\int_0^{\rho_M} \rho^{1-n}\,H^2(\rho)\, d\rho & \leq & \frac{1}{\Lambda_*}
\left ( \int_0^{\rho_M} \rho^{\frac{1-n}{p-1}} \left ( \frac{dH}{d\rho}
\right )^{\frac{p}{p-1}}\, d\rho \right )^{\frac{2(p-1)}{p}} \\ 
& = & \frac{1}{\Lambda_*} \left ( (n\omega_n)^{\frac{1}{p-1}}
\int_D \phi^p d\mu \right )^{\frac{2(p-1)}{p}} \\ 
& = & \frac{1}{\Lambda_*} (n \omega_n)^{2/p} \left ( 
\int_D \phi^p d\mu \right )^{\frac{2(p-1)}{p}},
\end {eqnarray*} 
with equality if and only if $D$ is a ball and $\phi$ is 
radial. Substituting this last inequality into \eqref{comparison-ode3}, 
we have 
\begin {eqnarray*} 
\left ( \int_D \phi^{p-1} d\mu \right )^2 & \leq & \frac{2n^2 \omega_n^{n/2}}
{p \calCp(D)} |D|^{\frac{n-2}{n}} \left ( \int_D \phi^p d\mu \right )
^{\frac{2(p-1)}{p}} \\ 
&&\qquad- (n-2) \omega_n^{\frac{2-n}{n}} |D|^{\frac{n-2}{n}}\frac{1}{\Lambda_*}
(n\omega_n)^{2/p} \left ( \int_D \phi^p d\mu \right )^{\frac{2(p-1)}{p}}.
\end {eqnarray*} 
Since $\Lambda_* = (n \omega_n)^{\frac{2-p}{p}} \calCp(D^*)$ by 
\eqref{comparison-ode7} and \eqref{lambda-star2}, the main inequality \eqref{main-ineq}
follows with equality if and only if $D$ is a ball.
\hfill $\square$


\end{document}